# Optimal Purchasing Policy For Mean-Reverting Items in a Finite Horizon

Alon Dourban and Liron Yedidsion
Technion – Israel Institute of Technology

**Abstract**

In this research we study a finite horizon optimal purchasing problem for items with a mean reverting price process. Under this model a fixed amount of identical items are bought under a given deadline, with the objective of minimizing the cost of their purchasing price and associated holding cost. We prove that the optimal policy for minimizing the expected cost is in the form of a time-variant threshold function that defines the price region in which a purchasing decision is optimal. We construct the threshold function with a simple algorithm that is based on a dynamic programming procedure that calculates the cost function. As part of this procedure we also introduce explicit equations for the crossing time probability and the overshoot expectation of the price process with respect to the threshold function. The characteristics and dynamics of the threshold function are analyzed with respect to time, holding cost, and different parameters of the price process, and yields meaningful practical insights, as well as theoretical insights.

## 1 Introduction

In this paper, we study an inventory management problem of items with a mean-reverting price process. We find the optimal policy for managing purchasing decisions when the underlying item's price is stochastic and tends to revert back to a long term average price. We analyze this system under a general convex holding cost and a finite deadline constraint. Our work introduces a simple purchasing policy that is based on a threshold value as a function of time. The threshold value indicates whether to purchase the item or to wait for a better price in the future. The problem we study can be found in the core operations of many businesses. For example, the model can be applied to manage an airport's fuel inventory. In this operation, the airport's flights schedule determines the deadline for the fuel demand, where the fuel is bought under a fluctuating price. Hence, the fuel inventory purchasing decision needs to be managed optimally in order to minimize the associated total costs that are composed of the direct price of the fuel and the holding costs. Our motivation to focus on items with stochastic prices relies on the growing need of businesses to take into account in their logistic operations the random behavior of their inventory purchasing price. Moreover, many evidence of mean-reverting behavior are detected in the random fluctuations of prices of popular goods, mainly commodities, see for example [4],[5],[26],[27],[28],[33]. Although our model primarily aims for optimal policy of classic logistic systems, it can be also adopted to the optimization of alternative systems that encounter mean-reversion in their underlying process. Such application include electricity markets [17], exchange rates [12], and pairs trading [9].

Researches from the optimal inventory control stream addressed the stochastic behavior of prices in several works. Kalymon [21] introduced a single item, multi period, inventory system for a system with random item price that follows a discrete Markovian process, and constructed a $(s(p); S(p))$ optimal policy. Golabi [14] considered a similar problem under an independent continuous price distribution framework. However non of these works offer closed form solutions to calculate the optimal policy in case of mean reverting price process. Berling and Martinez-de-Albeniz [3] work considers a closely related model to our research as they apply a numerical procedure in a continuous review inventory where the price follows the mean reverting Ornstein-Uhlenbeck process ($OU$). Goel and Gutierrez [13] also consider the $OU$ price process, and present a computational study that offers approximation for optimal inventory policy under multiechelon model in which purchasing can be made from the spot and forward markets. Other papers, that deal with stochastic prices in inventory systems include Gavirneni [11] and Wang [31] that study models where the price fluctuates myopically, and the work of Li and Kouvelis [22] which studies risk sharing contracts under a random price that follows the Geometric Brownian Motion ($GBM$) process, using a binomial approximation. For a comprehensive review of other related works that consider random price in inventory optimization see Haksöz and Seshadri [16].

The analysis of stochastic price process in optimal inventory control shares similar challenges to those found in optimal stopping problems. Optimal stopping methods are widely used in financial literature related to financial derivatives. Jacka [20] shows that the problem of American put option is equivalent to an optimal stopping problem. Since the seminal work of Black and Scholes [6], there is a rich research stream that has been developed in the financial literature regarding pricing and modeling of American options, See [18] for a review. However, in most cases there are no closed form solutions for this kind of problems, and therefore, different numerical procedures were developed. The majority of the financial research focus on an underlying process that follows the $GBM$ process, but some of the procedures that were developed can also capture mean reverting processes. This includes: the censored binomial tree of Nelson and Ramaswamy [24], the two dimensional tree of Hahn and Dyer [15], the trinomial tree of Hull and White [19], and the least squares simulation techniques of Logstaff and Schwartz [23].

Aiming to find the optimal policy of mean reverting process, some papers address sub-problems that are studied in our work as well. This includes the probability of the process crossing some given value, and the expected value of the process at this crossing time. Finster [10] and Novikov [25] studied asymptotic results for the first crossing time of a constant level for the first order autoregressive ($AR(1)$) process, and Alili et al. [1] presented three approximation methods for the crossing time of an $OU$ process. Christensen and Novikov [8] find conditions where the optimal stopping of a discounted $AR(1)$ process in an infinite horizon is in the form of a constant threshold function, and in a following work Christensen [7] assumed exponential random terms in the process, and find explicit distributions for the crossing time and expected value of the process under this assumption.

In this work, we introduce a simple optimal policy for the optimal inventory control of items with mean reverting price process. Our results give a practical tool for managing inventory in stochastic price systems, and in addition have important contribution to the theoretical study of optimal stopping of mean-reverting process. The optimal policy is based on an increasing threshold function in time, that its calculation is based on the crossing time probability of the process, and its expected value. We explicitly define equations for these terms, and find the value function of the system at a given time, using Bellman equation (see [2]). We prove the existence and basic properties of the threshold under this setting. The threshold function is calculated in a simple binary search that finds the indifference price between purchasing and waiting decisions. Lastly, we analyze the threshold dynamics with respect to logistic parameters of the inventory system, and to the price process parameters.

The remainder of the paper is organized as follows: The model is characterized in Section 2. The calculation of the crossing time probability and the overshoot expectation is constructed in Section 3, and the threshold dynamics analysis is presented in Section 4.

## 2 Model

In this research, we study a purchasing problem of items with mean reverting price process. Without loss of generality we consider a demand for a single item. The purchasing decision can be called in a finite time horizon of length $T$. Lead time is zero, and no backlogs are allowed. The item incurs a holding cost per time range $[t, s]$ denoted by, $h_{t,s}$. Our goal is to minimize the expected total cost of supplying the item within its planning horizon. The item's price, denoted by $X_t$, tends to revert back to a long run average price, $\theta$, according to the following autoregressive model of order 1 ($AR(1)$):

$$X_t = \theta \cdot K \cdot dt + (1 - K \cdot dt) \cdot X_{t-dt} + \sigma \cdot \varepsilon_t, \tag{1}$$

where $K$ is the reverting rate, $dt$ is the time interval length, and $\sigma$ is the degree of volatility. The stochastic term of the process is denoted by a series of i.i.d random variables with zero mean and variance $dt$. That is, $\varepsilon_t \sim N(0, dt)$. The parameters $K$ and $dt$ are constrained to hold: $|1 - dt \cdot K| < 1$ for keeping the price process stationary. This is in fact the discretized version of the arithmetic $OU$ [30] process, where the change in the process is described by the following equation:

$$dX_t = K\left(\theta - X_{t-dt}\right) dt + \sigma \cdot \varepsilon_t, \qquad (2)$$

In order to simplify notation to natural numbers indexes we also define: $n = t/dt$, and $t_n = n \cdot dt$ for $0 \leq n \leq N$, where $N = T/dt$. Under this notation we denote the price process as: $X_{t_0}, X_{t_1}, X_{t_2}, \ldots, X_{t_N}$, and we denote by $X_{t_n, t_i}(x)$, the random value of the process, $X_{t_i}$, conditioned that $X_{t_n} = x$, for $t_n < t_i$. According to eq. (1), we get:

$$\begin{aligned}
X_{t_{i-1}, t_i}(x) &= \theta \cdot K \cdot dt + (1 - K \cdot dt) \cdot x + \sigma \cdot \varepsilon_{t_i} \\
X_{t_{i-2}, t_i}(x) &= \theta \cdot K \cdot dt \cdot (1 + (1 - K \cdot dt)) + (1 - K \cdot dt)^2 \cdot x \\
&\quad + \sigma \cdot \left(\varepsilon_{t_i} + \varepsilon_{t_{i-1}} \cdot (1 - K \cdot dt)\right) \\
&= \theta \cdot dt \cdot K \cdot \sum_{j=0}^{1} (1 - dt \cdot K)^j + (1 - dt \cdot K)^2 \cdot x \\
&\quad + \sigma \cdot \sum_{j=0}^{1} \varepsilon_{t_{i-j}} \cdot (1 - dt \cdot K)^j,
\end{aligned}$$

which can be extended recursively to $X_{t_n, t_i}(x)$, for $0 \leq t_n < t_i \leq t_N$. This gives us:

$$\begin{aligned}
X_{t_n, t_i}(x) &= \theta \cdot dt \cdot K \cdot \sum_{j=0}^{i-n-1} (1 - dt \cdot K)^j + (1 - dt \cdot K)^{i-n} \cdot x \\
&\quad + \sigma \cdot \sum_{j=0}^{i-n-1} \varepsilon_{t_{i-j}} \cdot (1 - dt \cdot K)^j,
\end{aligned}$$

and by aggregating the geometric progression term, $\sum_{j=0}^{i-n-1} (1 - dt \cdot K)^j$, we get:

$$\begin{aligned}
X_{t_n, t_i}(x) &= \theta \cdot dt \cdot K \cdot \frac{\left(1 - (1 - dt \cdot K)^{i-n}\right)}{1 - (1 - dt \cdot K)} + (1 - dt \cdot K)^{i-n} \cdot x \\
&\quad + \sigma \cdot \sum_{j=0}^{i-n-1} \varepsilon_{t_{i-j}} \cdot (1 - dt \cdot K)^j \\
&= \theta + (x - \theta)(1 - dt \cdot K)^{i-n} + \sigma \sum_{j=0}^{i-n-1} (1 - dt \cdot K)^j \varepsilon_{t_{i-j}}. \qquad (3)
\end{aligned}$$

Note that we do not restrict the price process, $X_t$, to obtain only positive values, as it may represent the spread of two related assets in a pairs trading model (Ekstrom et al. 2011). In addition, we consider a quite general time dependent holding cost. Our only restriction is that the holding cost is convex in its time range. However, in our model the only relevant holding cost terms are $h_{t_n, t_N}$, for $0 \leq t_n \leq t_{N-1}$. Therefore, any holding cost function can be degenerated into an additive function of the set of holding cost per time unit, denoted $h_{t_n}$, where $h_{t_n}$ decreases (weakly) in $t_n$.

We show in Section 3 that the optimal policy is based on a threshold function, in which at every time point, $t$, there exists a single threshold value that states the following purchasing policy: buy if the item's price drops below the threshold, otherwise wait. In case the item was not purchased until time $T$, it is bought at time $T$ regardless of the value of $X_T$. Accordingly, a purchasing policy $\pi$ defines a set of threshold values, $\mathbf{b}_\pi = \{b_\pi(t_0), \ldots, b_\pi(t_N)\}$. Let $V_\pi(x, t)$ denote the value function of the system at time $t$, and price $x$, under policy $\pi$. That is, $V_\pi(x, t)$ defines the expected future total cost of the system. Thus, we aim to find an optimal policy, denoted as $\pi^*$, which minimizes the expected total cost of the system at any time within the time horizon, given the market price at that time. This optimal policy is denoted by $V(x,t)$. That is,

$$\pi^* = \arg\min V_\pi(x,t),$$

where we define the optimal terms as:

$$b(t) \equiv \pi^*,$$
$$V(x,t) \equiv V_{\pi^*}(x,t).$$

We define the set of prices that are higher than the threshold value at time, $t$, as the *continuation region* at time $t$, and the set of prices that are lower than the threshold value at time, $t$, as the *stopping region* at time $t$.

## 3 Dynamic Programming based solution

Our aim is to find the optimal policy for the purchasing problem at any point in the time horizon $0 \leq t \leq T$. A natural way to approach this problem is by applying a Dynamic Programming based solution according to Bellman's equation, where the purchasing constraint at time $T$ acts as the stopping condition. That is,

$$V(x,T) = x, \tag{4}$$

as at time $T$ there is no other choice rather than to buy the item at price $x$. Eq. (4) is used to calculate the value function at earlier time points, $t_0 \leq t_n < t_N$ recursively:

$$V(x,t_n) = \min\left\{E\left[V(X_{t_n,t_{n+1}}(x),t_{n+1})\right], x + dt \cdot \sum_{i=n}^{N-1} h_{t_i}\right\}. \tag{5}$$

The left term inside the minimum function, $E\left[V(X_{t_n,t_{n+1}}(x),t_{n+1})\right]$, represents the expected total cost of the system when we wait with the purchasing decision at time $t$, and continue according to the optimal policy in the reminder of the horizon. We define this expected value as the value function under *continuation* decision, noted as $V^C(x,t_n)$. The right term inside the minimum function, $x + dt \cdot \sum_{i=n}^{N-1} h_{t_i}$, represents the total cost of the system when the item is purchased at time $t_n$ at price $x$, and a holding cost of $dt \cdot \sum_{i=n}^{N} h_i$ is added to the item's price. At each time interval, the price region is separated into *continuation region* in which $V^C(x,t_n) < x + dt \cdot \sum_{i=n}^{N-1} h_{t_i}$, and a *stopping region* in which $V^C(x,t_n) > x + dt \cdot \sum_{i=n}^{N-1} h_{t_i}$. The threshold value, $b(t_n)$, represents the indifference price between these two decisions, defined by:

$$V^C(b(t_n),t_n) = b(t_n) + dt \cdot \sum_{i=n}^{N-1} h_{t_i} = V(b(t_n),t_n). \tag{6}$$

Note that according to this terminology the threshold price at time $t_N$ does not exist, as there is no defined value from continuation, but for simplicity of the calculations we define $b(t_N) = \infty$.

In the last period before the end of the time horizon, $t_{N-1}$, the decision maker has to decide whether to purchase the item at the price of $x$ and pay $dt \cdot h_{t_{N-1}}$ holding cost, or to wait with the purchasing decision to the end of the time horizon where the item would be purchased at the price of $X_{t_N}$ with no additional holding cost. Hence, according to eqs. (2) and (5), $b(t_{N-1})$ satisfies

$$E\left[b(t_{N-1}) + K(\theta - b(t_{N-1}))dt + \sigma \cdot \varepsilon_t\right] = b(t_{N-1}) + dt \cdot h_{t_{N-1}}, \tag{7}$$

and as $E\left[\varepsilon_t\right] = 0$, eq. (7) gives

$$b\left(t_{N-1}\right) + K\left(\theta - b\left(t_{N-1}\right)\right) dt = b\left(t_{N-1}\right) + dt \cdot h_{t_{N-1}}.$$

Thus, we get the following Corollary:

**Corollary 1**

$$b\left(t_{N-1}\right) = \theta - h_{t_{N-1}}/K.$$

The threshold value one time interval before the deadline is given in a closed form in Corollary 1. However, the calculation procedure of threshold values in earlier times is more complex. First, we have to show that there is a single threshold at each time period. To do so, we show that the following Properties hold:

1. At any time point, $0 \leq t_n \leq t_{N-1}$, there exists a finite price, $x^H(t_n) \equiv \theta - h_{t_{N-1}}/K$, for which any $x \geq x^H(t_n)$ satisfies $x + dt \cdot \sum_{i=n}^{N-1} h_{t_i} > V^C(x, t_n)$.

2. At any time point, $0 \leq t_n \leq t_{N-1}$, there exists a finite price, $x^L(t_n) \equiv \min\left\{\frac{b(t_{n+1}) - \theta dt \cdot K}{(1 - dt \cdot K)}, \frac{\theta \cdot dt \cdot K - \sigma \cdot \sqrt{dt} - h_{t_n} \cdot dt}{dt \cdot K}\right\}$, for which any $x \leq x^L(t_n)$ satisfies $x + dt \cdot \sum_{i=n}^{N} h_i < V^C(x, t_n)$.

3. $V^C(x, t_n)$ is concave and increasing in $x$.

Then, we rely on properties $1 - 3$, and derive that at any time point there exists a single threshold value, $b(t_n)$. This is proved in the following Lemma:

**Lemma 1** *Properties* $1 - 3$ *hold, and there exists a finite single threshold value, $b(t_n)$, for $0 \leq t_n < t_N < \infty$.*

**Proof.** See proof in Appendix 1. ∎

Next, we define preliminary terms that are used to calculate the value function: the *Crossing Time Probability* and *Overshoot Expectation*.

**Definition 1** *Crossing Time*

We define the earliest time point, after time $t_n$, in which the price process is lower than its respective threshold, as the Crossing Time of the price process, denoted by $\tau_{t_n}(x)$. That is,
$\tau_{t_n}(x) = \inf\{t_i > t_n : X_{t_n, t_i}(x) \leq b(t_i)\}$, for $0 \leq t_n < t_i \leq t_N$.

**Definition 2** *Crossing Time Probability*

We define the Crossing Time Probability at time $t_i$: the probability that $\tau_{t_n}(x)$ occurs at time $t_i$. We denote the Crossing Time Probability by $P_{t_n, t_i}(x)$. That is,

$$P_{t_n, t_i}(x) = \Pr\left(\tau_{t_n}(x) = t_i\right), \text{ for } 0 \leq t_n < t_i \leq t_N.$$

**Definition 3** *Overshoot Expectation*

We define the Overshoot Expectation at time $t_i$: the expected value of the process at $\tau_{t_n}(x)$ conditioned that the Crossing Time occurs at time $t_i$ and that the price at time $t_n$ is $x$, for $0 \leq t_n < t_i \leq t_N$. We denote the Overshoot Expectation by $E_{t_n, t_i}(x)$. That is,

$$E_{t_n, t_i}(x) = E\left(X_{t_n, t_i}(x) \mid \tau_{t_n}(x) = t_i\right).$$

By using $P_{t_n, t_i}(x)$ and $E_{t_n, t_i}(x)$ for $0 \leq t_n < t_i \leq t_N$ we can construct the value function in the *continuation region*, for time $t_n$, and price $x > b(t_n)$ according to the following recursive equation:

$$\begin{aligned} V(x, t_n) = V^C(x, t_n) &= P_{t_n, t_{n+1}}(x) \cdot \left(E_{t_n, t_{n+1}}(x) + dt \cdot \sum_{j=n+1}^{N-1} h_{t_i}\right) \\ &+ \left(1 - P_{t_n, t_{n+1}}(x)\right) \cdot E\left[V(X_{t_n, t_{n+1}}(x), t_{n+1}) \mid X_{t_n, t_{n+1}}(x) > b(t_{n+1})\right] \\ &= \sum_{i=n+1}^{N} P_{t_n, t_i}(x) \cdot \left(E_{t_n, t_i}(x) + dt \cdot \sum_{j=i}^{N-1} h_{t_i}\right), \end{aligned} \qquad (8)$$

and in the *stopping region* for time $t_n$ and price $x \leq b(t_n)$ we get:

$$V(x, t_n) = x + dt \cdot \sum_{i=n}^{N-1} h_{t_i}.$$

## 3.1 Crossing Time Probability and Overshoot Expectation calculation

After defining the terms for *Crossing Time Probability* and *Overshoot Expectation* in definitions 2 and 3, in this section we present their calculation. The calculation of $P_{t_n,t_i}(x)$ and $E_{t_n,t_i}(x)$ is based on decomposing these terms to explicit functions of multivariate normal variables, and then calculate their probability function, and truncated expectation.

According to eq. (3), $X_{t_n,t_i}(x)$, is in fact a random variable with normal distribution. That is,

$$X_{t_n,t_i}(x) \sim N\left(\mu_{t_i,t_n}(x), \sigma^2_{t_i,t_n}\right)$$

which according to eq. (3) has expectation, $\mu_{t_n,t_i}(x)$, that equals:

$$\mu_{t_n,t_i}(x) = \theta + (x - \theta)(1 - dt \cdot K)^{i-n}, \tag{9}$$

and variance, $\sigma^2_{t_n,t_i}$, that equals:

$$\sigma^2_{t_n,t_i} = Var\left[\sigma \sum_{j=0}^{i-n-1} (1 - dt \cdot K)^j \varepsilon_{t_{i-j}}\right]. \tag{10}$$

$\varepsilon_{t_{n+1}}, \varepsilon_{t_{n+2}}, \ldots, \varepsilon_{t_i}$ are i.i.d random variables with variance that equals $dt$. Thus, by aggregating the geometric progression in eq. (10), we get:

$$\sigma^2_{t_n,t_i} = dt \cdot \sigma^2 \frac{1 - (1 - dt \cdot K)^{2(i-n)}}{1 - (1 - dt \cdot K)^2}. \tag{11}$$

In correspondence to $X_{t_n,t_i}(x)$, we define its relative standard normal variable by

$$Z_{t_n,t_i}(x) := \frac{X_{t_n,t_i}(x) - \mu_{t_n,t_i}(x)}{\sigma_{t_n,t_i}}, \tag{12}$$

where $Z_{t_n,t_i}(x) \sim N(0,1)$. The set of random variables, $Z_{t_n,t_{n+1}}(x), \ldots, Z_{t_n,t_i}(x)$, are correlated with a symmetric covariance matrix denoted by $\Sigma_{t_{n+1},t_i}$. Each covariance element, $cov(Z_{t_n,t_l}(x), Z_{t_n,t_k}(x))$, for $t_{n+1} \leq t_l \leq t_k \leq t_i$, in the covariance matrix, $\Sigma_{t_{n+1},t_i}$, represents the covariance between $Z_{t_n,t_l}(x)$ and $Z_{t_n,t_k}(x)$, and equals:

$$cov(Z_{t_n,t_l}(x), Z_{t_n,t_k}(x)) = E[Z_{t_n,t_l}(x) \cdot Z_{t_n,t_k}(x)] - E[Z_{t_n,t_l}(x)] \cdot E[Z_{t_n,t_k}(x)],$$

where $E[Z_{t_n,t_l}(x)] \cdot E[Z_{t_n,t_k}(x)] = 0$. According to eqs. (12), and (9):

$$cov(Z_{t_n,t_l}(x), Z_{t_n,t_k}(x)) = \frac{1}{\sigma_{t_n,t_l} \cdot \sigma_{t_n,t_k}} \cdot E\left[\left(\sigma \sum_{j=0}^{l-n-1}(1-dt\cdot K)^j \varepsilon_{t_{l-j}}\right) \cdot \left(\sigma \sum_{r=0}^{k-n-1}(1-dt\cdot K)^r \varepsilon_{t_{k-r}}\right)\right]. \tag{13}$$

Note that $E[\varepsilon_{t_{l-j}} \cdot \varepsilon_{t_{k-r}}] = dt$ for $t_{l-j} = t_{k-r}$, and $E[\varepsilon_{t_{l-j}} \cdot \varepsilon_{t_{k-r}}] = 0$ for $t_{l-j} \neq t_{k-r}$. $t_{l-j} = t_{k-r}$ for pairs $(j,r) \in \{(0, k-l), \ldots, (l-n-1, k-n-1)\}$. Therefore,

$$cov(Z_{t_n,t_l}(x), Z_{t_n,t_k}(x)) = \frac{dt \cdot \sigma^2}{\sigma_{t_n,t_l} \cdot \sigma_{t_n,t_k}} \cdot \sum_{j=0}^{l-n-1} (1 - dt \cdot K)^{k-l+2\cdot j}.$$

By aggregating the geometric progression we get:

$$cov\left(Z_{t_n,t_l}(x), Z_{t_n,t_k}(x)\right) = \frac{dt \cdot \sigma^2 \cdot (1 - dt \cdot K)^{k-l}}{\sigma_{t_n,t_l} \cdot \sigma_{t_n,t_k}} \cdot \frac{1 - (1 - dt \cdot K)^{2(l-n)}}{1 - (1 - dt \cdot K)^2}, \text{ for } t_{n+1} \leq t_l \leq t_k \leq t_i. \quad (14)$$

Next, we define the respective standardized threshold values. That is, the corresponding values of $b(t_n)$ that are adjusted to the standardized process $Z_{t_n,t_i}(x)$. We define

$$\mathbf{B}_{t_n,t_i}(x) := \left[\beta_{t_n,t_{n+1}}(x), \beta_{t_n,t_{n+2}}(x), ..., \beta_{t_n,t_i}(x)\right],$$

as the vector of standardized threshold values for time points $t_{n+1}, \ldots, t_i$, where each element, $\beta_{t_n,t_{n+l}}(x)$, for $t_{n+1} \leq t_{n+l} \leq t_i$, is set to the corresponding standardized value of $b(t_{n+l})$ with respect to the standardized process $Z_{t_n,t_{n+l}}(x)$. That is,

$$\beta_{t_n,t_l}(x) := \frac{b(t_l) - \left[\theta + (x - \theta) \cdot (1 - dtK)^{l-n}\right]}{\sigma_{t_n,t_l}}, \text{ for } n+1 \leq l \leq i. \quad (15)$$

Similarly, we define $\widehat{\mathbf{B}}_{t_n,t_i}(x)$ as the vector of standardized threshold values for time points $t_{n+1}, \ldots, t_i$, with the exception that the element that corresponds to time $t_i$ is set to $-\infty$. That is,

$$\widehat{\mathbf{B}}_{t_n,t_i}(x) := \left[\widehat{\beta}_{t_n,t_{n+1}}(x), \widehat{\beta}_{t_n,t_{n+2}}(x), ..., \widehat{\beta}_{t_n,t_i}(x)\right],$$

where

$$\widehat{\beta}_{t_n,t_l}(x) := \begin{cases} \beta_{t_n,t_l}(x), & \text{for } n+1 \leq l \leq i-1, \\ -\infty, & l = i \end{cases}.$$

**Lemma 2**

$$P_{t_n,t_i}(x) = F\left(-\widehat{\mathbf{B}}_{t_n,t_i}(x), \Sigma_{t_n,t_i}\right) - F\left(-\mathbf{B}_{t_n,t_i}(x), \Sigma_{t_n,t_i}\right), \quad (16)$$

where $F(\mathbf{a}^n, \Sigma^{n \times n})$ is the cumulative distribution function of a standard multivariate normal variable of order $n$, on the vector $\mathbf{a} \in R^n$, with an $n \times n$ covariance matrix $\Sigma$.

**Proof.** According to Definition 2, $P_{t_n,t_i}(x)$ can be formulated as the probability that the process is above the threshold at times $\{t_{n+1}, \ldots, t_{i-1}\}$, and below the threshold at $t_i$. That is,

$$P_{t_n,t_i}(x) = \Pr\left(X_{t_n,t_i}(x) \leq b(t_i), X_{t_n,t_{i-1}}(x) > b(t_{i-1}), ..., X_{t_n,t_{n+1}}(x) > b(t_{n+1})\right).$$

As $\beta_{t_n,t_{n+1}}(x), ..., \beta_{t_n,t_i}(x)$ denote the corresponding standardized values of $b(t_{n+1}), \ldots, b(t_i)$, we get that,

$$P_{t_n,t_i}(x) = \Pr\left(Z_{t_n,t_i}(x) \leq \beta_{t_n,t_i}(x), Z_{t_n,t_{i-1}}(x) > \beta_{t_n,t_{i-1}}(x), ..., Z_{t_n,t_{n+1}}(x) > \beta_{t_n,t_{n+1}}(x)\right).$$

By the symmetry of $Z_{t_n,t_{n+1}}(x), \ldots, Z_{t_n,t_i}(x)$ and the law of total probability we get:

$$P_{t_n,t_i}(x) = \Pr\left(Z_{t_n,t_{i-1}}(x) \leq -\beta_{t_n,t_{i-1}}(x), ..., Z_{t_n,t_{n+1}}(x) \leq -\beta_{t_n,t_{n+1}}(x)\right)$$
$$- \Pr\left(Z_{t_n,t_i}(x) \leq -\beta_{t_n,t_i}(x), Z_{t_n,t_{i-1}}(x) \leq -\beta_{t_n,t_{i-1}}(x), ..., Z_{t_n,t_{n+1}}(x) \leq -\beta_{t_n,t_{n+1}}(x)\right)$$
$$= F\left(-\widehat{\mathbf{B}}_{t_n,t_i}(x), \Sigma_{t_n,t_i}\right) - F\left(-\mathbf{B}_{t_n,t_i}(x), \Sigma_{t_n,t_i}\right).$$

∎

Next, we need to calculate the *Overshoot Expectation* of the process. This is done by forming $E_{t_n,t_i}(x)$ in terms of the conditional expectation function of a multivariate normal variable. First, we define the vectors:

$$\mathbf{A}_{t_n,t_i.t_j}(x) = \left[\alpha_{t_n,t_{n+1}.t_j}(x), ..., \alpha_{t_n,t_{j-1}.t_j}(x), \alpha_{t_n,t_{j+1}.t_j}(x), ..., \alpha_{t_n,t_i.t_j}(x)\right], \quad (17)$$

$$\widehat{\mathbf{A}}_{t_n,t_i.t_j}(x) = \left[\widehat{\alpha}_{t_n,t_{n+1}.t_j}(x), ..., \widehat{\alpha}_{t_n,t_{j-1}.t_j}(x), \widehat{\alpha}_{t_n,t_{j+1}.t_j}(x), ..., \widehat{\alpha}_{t_n,t_i.t_j}(x)\right], \quad (18)$$

where

$$\alpha_{t_n,t_l.t_j}(x) = \frac{\left(\beta_{t_n,t_l}(x) - Cov\left(Z_{t_n,t_l}(x), Z_{t_n,t_j}(x)\right) \cdot \beta_{t_n,t_j}(x)\right)}{\sqrt{1 - Cov\left(Z_{t_n,t_l}(x), Z_{t_n,t_j}(x)\right)^2}}, \quad (19)$$

$$\text{for } n+1 \leq l \leq i, \text{ and } n+1 \leq j \leq i, l \neq j,$$

and

$$\widehat{\alpha}_{t_n,t_l.t_j}(x_n) = \begin{cases} \alpha_{t_n,t_l.t_j}(x_n), & \text{for } n+1 \leq l \leq i-1, \text{ and } n+1 \leq j \leq i, l \neq j, \\ -\infty, & \text{for } l = i, \text{ and } n+1 \leq j \leq i, i \neq j. \end{cases} \quad (20)$$

We define the matrix, $\mathbf{M}_{t_i,t_n.t_j}(x)$, as the $(i-n-1) \times (i-n-1)$ first-order partial correlation matrix of $\{Z_{t_n,t_{n+1}}(x), \ldots, Z_{t_n,t_i}(x)\}$, for $n+1 \leq j \leq i$, $i \neq j$, when removing the controlling variable $Z_{t_n,t_j}(x)$. That is, the $\{t_l, t_m\}$ entry to $\mathbf{M}_{t_i,t_n.t_j}(x)$, denoted by $\rho_{t_l,t_m.t_j}(x)$, for $n+1 \leq l \leq i$, $n+1 \leq m \leq i$, $n+1 \leq j \leq i$, $l,m \neq j$, equals:

$$\rho_{t_l,t_m.t_j}(x) = \frac{cov\left(Z_{t_n,t_l}(x), Z_{t_n,t_m}(x)\right) - cov\left(Z_{t_n,t_l}(x), Z_{t_n,t_j}(x)\right) \cdot cov\left(Z_{t_n,t_m}(x), Z_{t_i,t_m}(x)\right)}{\sqrt{1 - cov\left(Z_{t_n,t_l}(x), Z_{t_n,t_j}(x)\right)^2} \cdot \sqrt{1 - cov\left(Z_{t_n,t_m}(x), Z_{t_n,t_j}(x)\right)^2}}. \quad (21)$$

With these terms defined we use Tallis' equation (see [29]), for the mean of a standardized normal truncated multivariate distribution, in order to calculate the expected value of the process at time $t_n$, conditioned that $X_{t_n} = x$, and that the corresponding *Crossing Time* is triggered after time, $t_j$, for $0 \leq t_n \leq t_j \leq t_i$, denoted as $E\left(X_{t_n,t_i}(x) | \tau_{t_n}(x) > t_j\right)$. We obtain this term for $t_j = t_i$, $t_j = t_{i-1}$, and use them to derive $E_{t_n,t_i}(x)$. Applying Tallis' equation with the set of $\{Z_{t_n,t_{n+1}}(x), \ldots, Z_{t_n,t_i}(x)\}$ as the standardized normal multivariate variable truncated respectively on the coordinates $\{\beta_{t_n,t_{n+1}}(x_n), \ldots, \beta_{t_n,t_{n_i}}(x_n)\}$ gives:

$$E\left(Z_{t_n,t_i}(x) | \tau_{t_n}(x) > t_i\right) \quad (22)$$
$$= \frac{\sum_{l=n+1}^{i} cov\left(Z_{t_n,t_i}(x), Z_{t_n,t_l}(x)\right) \cdot \phi\left(\beta_{t_n,t_l}(x_n)\right) \cdot F\left(-\mathbf{A}_{t_n,t_i.t_l}(x_n), \mathbf{M}_{t_i,t_n.t_l}(x)\right)}{F\left(-\mathbf{B}_{t_n,t_i}(x_n), \Sigma_{t_n,t_i}\right)},$$

$$E\left(Z_{t_n,t_i}(x) | \tau_{t_n}(x) > t_{i-1}\right) \quad (23)$$
$$= \frac{\sum_{l=n+1}^{i} cov\left(Z_{t_n,t_i}(x), Z_{t_n,t_l}(x)\right) \cdot \phi\left(\widehat{\beta}_{t_n,t_l}(x_n)\right) \cdot F\left(-\widehat{\mathbf{A}}_{t_n,t_i.t_l}(x_n), \mathbf{M}_{t_i,t_n.t_l}(x)\right)}{F\left(-\widehat{\mathbf{B}}_{t_n,t_i}(x_n), \Sigma_{t_n,t_i}\right)},$$

where $\phi(x)$ is the probability density function ($PDF$) of a standard normal variable.

Taking the standardized terms, $E\left(Z_{t_n,t_i}(x) | \tau_{t_n}(x) > t_i\right)$, $E\left(Z_{t_n,t_i}(x) | \tau_{t_n}(x) > t_{i-1}\right)$ and transforming them to the form of $E\left(X_{t_n,t_i}(x) | \tau_{t_n}(x) > t_i\right)$, $E\left(X_{t_n,t_i}(x) | \tau_{t_n}(x) > t_{i-1}\right)$ gives:

$$E\left(X_{t_n,t_i}(x) | \tau_{t_n}(x) > t_i\right) \quad (24)$$
$$= \mu_{t_n,t_i}(x) + \frac{\sigma_{t_n,t_i} \cdot \sum_{l=n+1}^{i} cov\left(Z_{t_n,t_i}(x), Z_{t_n,t_l}(x)\right) \cdot \phi\left(\beta_{t_n,t_l}(x_n)\right) \cdot F\left(-\mathbf{A}_{t_n,t_i.t_l}(x_n), \mathbf{M}_{t_i,t_n.t_l}(x)\right)}{F\left(-\mathbf{B}_{t_n,t_i}(x_n), \Sigma_{t_n,t_i}\right)},$$

$$E\left(X_{t_n,t_i}(x) | \tau_{t_n}(x) > t_{i-1}\right) \quad (25)$$
$$= \mu_{t_n,t_i}(x) + \frac{\sigma_{t_n,t_i} \cdot \sum_{l=n+1}^{i} cov\left(Z_{t_n,t_i}(x), Z_{t_n,t_l}(x)\right) \cdot \phi\left(\widehat{\beta}_{t_n,t_l}(x_n)\right) \cdot F\left(-\widehat{\mathbf{A}}_{t_n,t_i.t_l}(x_n), \mathbf{M}_{t_i,t_n.t_l}(x)\right)}{F\left(-\widehat{\mathbf{B}}_{t_n,t_i}(x_n), \Sigma_{t_n,t_i}\right)}.$$

$E_{t_n,t_i}(x)$ is derived by completing the total expectation equation when relying on the terms of

$E\left(X_{t_n,t_i}(x) \,|\, \tau_{t_n}(x) > t_i\right)$, $E\left(X_{t_n,t_i}(x) \,|\, \tau_{t_n}(x) > t_{i-1}\right)$ in eqs. (24),(25). In order to set the total expectation equation, we specify the conditional probability of the crossing time, denoted by $P_{t_n,t_i}\left(x \,|\, \tau_{t_n}(x) > t_{i-1}\right)$: the probability that $\tau_{t_n}(x)$ occurs at time $t_i$, conditioned that $\tau_{t_n}(x) > t_{i-1}$. According to Bayes rule, and eq. (16), we get that

$$P_{t_n,t_i}\left(x \,|\, \tau_{t_n}(x) > t_{i-1}\right) = \frac{P_{t_n,t_i}(x)}{F\left(-\widehat{\mathbf{B}}_{t_n,t_i}(x_n), \Sigma_{t_n,t_i}\right)}. \tag{26}$$

Lastly, we set the the total expectation equation, and and derive $E_{t_n,t_i}(x)$:

$$E_{t_n,t_i}(x) = \frac{\left[E\left(X_{t_n,t_i}(x) \,|\, \tau_{t_n}(x) > t_{i-1}\right) - \left(1 - P_{t_n,t_i}\left(x \,|\, \tau_{t_n}(x) > t_{i-1}\right)\right) \cdot E\left(X_{t_n,t_i}(x) \,|\, \tau_{t_n}(x) > t_i\right)\right]}{P_{t_n,t_i}\left(x \,|\, \tau_{t_n}(x) > t_{i-1}\right)}. \tag{27}$$

## 3.2 Threshold Function

In this section, we present the threshold calculation algorithm. But first, we prove that the threshold function maintains monotonicity in time, as stated in the following Lemma:

**Lemma 3** $b(t)$ *is monotonically increasing in* $t$, *for* $0 < t < T$.

**Proof.** Let $\pi_{-1}$ define the policy of setting the threshold at time $t_n$ to the optimal threshold at time $t_{n+1}$. That is, policy $\pi_{-1}$ can be defined as

$$\pi_{-1} \equiv b_{\pi_{-1}}(t_n) = b(t_{n+1}), \text{ for } 0 < t_n \leq t_{N-1}.$$

Let $V_{\pi_{-1}}(x, t_n)$ define the value function under policy $\pi_{-1}$. According to eq. (8) the value function under *continuation* decision under policy $\pi_{-1}$, denoted $V_{\pi_{-1}}^C(x, t_n)$, can be formulated as:

$$V_{\pi_{-1}}^C(x, t_n) = \sum_{i=n+1}^{N} P_{t_n,t_i}^{\pi_{-1}}(x) \cdot \left(E_{t_n,t_i}^{\pi_{-1}}(x) + dt \cdot \sum_{j=i}^{N-1} h_{t_j}\right),$$

where $P_{t_n,t_i}^{\pi_{-1}}(x)$ and $E_{t_n,t_i}^{\pi_{-1}}(x)$ stand for the *Crossing Time Probability* and the *Overshoot Expectation* under policy $\pi_{-1}$, respectively. Note that $P_{t_n,t_i}^{\pi_{-1}}(x) = P_{t_{n+1},t_{i+1}}(x)$, and $E_{t_n,t_i}^{\pi_{-1}}(x) = E_{t_{n+1},t_{i+1}}(x)$, but purchasing under policy $\pi_{-1}$ adds one additional time unit with holding cost compared to the optimal policy. Therefore, we get that

$$V_{\pi_{-1}}^C(x, t_n) = V^C(x, t_{n+1}) + \sum_{i=n+1}^{N-1} P_{t_{n+1},t_{i+1}}(x) \cdot dt \cdot h_{t_i}. \tag{28}$$

As our model restricts $h_t$ to decrease in $t$ we get from eq. (28):

$$V_{\pi_{-1}}^C(x, t_n) \leq V^C(x, t_{n+1}) + dt \cdot h_{t_n}, \tag{29}$$

and substituting for $V^C(x, t_{n+1})$ in eq. (6) into eq. (29) gives:

$$V_{\pi_{-1}}^C(b(t_{n+1}), t_n) \leq b(t_{n+1}) + dt \cdot \sum_{i=n+1}^{N-1} h_{t_i} + dt \cdot h_{t_n} = b(t_{n+1}) + dt \cdot \sum_{i=n}^{N-1} h_{t_i}. \tag{30}$$

Note that the *RHS* equals to the cost of a purchasing decision at time $t_n$ and price $b(t_{n+1})$. In addition, as $\pi_{-1}$ is a feasible policy, we get that $V_{\pi_{-1}}^C(b(t_{n+1}), t_n) \geq V^C(b(t_{n+1}), t_n)$. Combining this with eq. (30), we get that

$$V^C(b(t_{n+1}), t_n) \leq b(t_{n+1}) + dt \cdot \sum_{i=n}^{N-1} h_{t_i},$$

which together with the uniqueness of the threshold value (Lemma 1) implies that

$$b(t_n) \leq b(t_{n+1}), \text{ for } 0 < t < T.$$

∎

Note, that a sub-result of Lemma 3 is a new upper bound for the threshold function, which is tighter than the one that was suggested in Lemma 1. After establishing the *Crossing Time Probability* and *Overshoot Expectation* in eqs. (16) and (27), the threshold can be calculated using a simple recursive procedure. The value function in its continuation region, $V^C(x, t_n)$, is calculated according to eqs. (8), (16) and (27), and the threshold, $b(t_n)$, is determined by finding the value of $b(t_n)$ that satisfies eq. (6). However, the calculations of $P_{t_n,t_i}(x)$ and $E_{t_n,t_i}(x)$ for $t_{n+1} \leq t_i \leq t_N$, are based on the values of $\{b(t_{n+1}), \ldots, b(t_i)\}$. Therefore, we start calculating the threshold values from the end of the horizon, at time $t_N$, and move in sequential steps in length $dt$ backwards in time to $t_{n+1}$. As a stopping condition, at time $t_N$ the threshold, $b(t_N) = \infty$, and according to Corollary 1, $b(t_{N-1}) = \theta - h_{t_{N-1}}/K$. The reminder of the threshold values are determined in each step by a simple binary search method that finds the threshold value under a respective desired precision resolution, $\epsilon$. As the threshold is unique and holds upper and lower bounds, and is increasing in $t$ (Lemma 1, and Lemma 3), the binary search guarantees converging towards an optimal solution. This procedure is presented by the following algorithm:

**Algorithm 1** *Threshold Function*

1. Set $n = N - 2, b(t_N) = \infty, b(t_{N-1}) = \theta - h_{t_{N-1}}/K, LB = \frac{\theta \cdot dt \cdot K - \sigma \cdot \sqrt{dt} - h_{t_{n-2}} \cdot dt}{dt \cdot K}, UB = b(t_{N-1})$, and $x = (UB - LB)/2$

2. If $n \geq 0$

3. Calculate $V^C(x, t_n)$ according to eq. (8).

    3.1. while $\left| x + dt \cdot \sum_{i=n}^{N-1} h_{t_i} - V^C(x, t_n) \right| > \epsilon$

        3.1.1. if $x + dt \cdot \sum_{i=n}^{N-1} h_{t_i} > V^C(x, t_n)$, set $UB = x$
        3.1.2. else, set $LB = x$
        3.1.3. $x = (UB - LB)/2$
        3.1.4. Calculate $V^C(x, t_n)$ according to eq. (8).

    3.2. $b(t_n) = x$
    3.3. $n := n - 1$
    3.4. $LB = \min\left\{ \frac{b(t_{n+1}) - \theta dt \cdot K}{(1 - dt \cdot K)}, \frac{\theta \cdot dt \cdot K - \sigma \cdot \sqrt{dt} - h_{t_n} \cdot dt}{dt \cdot K} \right\}, UB = b(t_{n+1})$
    3.5. Go to step 2

4. end

## 4 Threshold Properties

In this section, we introduce some interesting properties of the threshold function. Within this section, we use the notation $b(t_n; \Lambda = \mathbf{c}), X_{t_n,t_i}(x; \Lambda = \mathbf{c})$ where $\Lambda$ is one or more of the model's parameters, and $\mathbf{c}$ corresponds to its value. This notation represents the threshold function and the price process respectively, at time $t_n$, under specific values, $\mathbf{c}$, of the noted parameters, $\Lambda$.

**Lemma 4** $b(t_n; \theta = \theta_0) + C = b(t_n; \theta = \theta_0 + C)$ for a constant $C$, and for $t_1 \leq t_n \leq t_N$.

**Proof.** According to the price process term in eq. (3), we get that $\theta$ shifts the price process only by a constant:

$$\begin{aligned} X_{t_n,t_i}(x+C;\theta=\theta_0+C) &= \theta_0 + C \\ &\quad + (x+C-(\theta_0+C))(1-dt\cdot K)^{i-n} \\ &\quad + \sigma \sum_{j=0}^{i-n-1} (1-dt\cdot K)^j \varepsilon_{t_{i-j}} \\ &= X_{t_n,t_i}(x;\theta=\theta_0) + C. \end{aligned} \quad (31)$$

Let $V^C(x,t_i;\theta=\theta_0), b(t_n;\theta=\theta_0)$ be the value function under *continuation* decision, and the threshold function for a process with $\theta=\theta_0$, respectively. Note that according to eqs. $(16),(24),(25),(27)$, and $(31)$ we get that for $0 \leq t_n < t_i \leq t_N$, and all $t_i < t \leq t_N$

$$E_{t_n,t_i}(x+C;\theta=\theta_0+C,\mathbf{b}=\mathbf{b}(;\theta=\theta_0)+C) = E_{t_n,t_i}(x;\theta=\theta_0,\mathbf{b}=\mathbf{b}(;\theta=\theta_0)) + C,$$

and

$$P_{t_n,t_i}(x+C;\theta=\theta_0+C,\mathbf{b}=\mathbf{b}(;\theta=\theta_0)+C) = P_{t_n,t_i}(x;\theta=\theta_0,\mathbf{b}=\mathbf{b}(;\theta=\theta_0)).$$

Hence, by eq. (8) we get that

$$V^C(x+C,t_n;\theta=\theta_0+C,\mathbf{b}=\mathbf{b}(;\theta=\theta_0)+C) = V^C(x,t_n;\theta=\theta_0,\mathbf{b}=\mathbf{b}(;\theta=\theta_0)) + C.$$

Therefore, as

$$V^C(b(t_n;\theta=\theta_0),t_n;\theta=\theta_0,\mathbf{b}=\mathbf{b}(;\theta=\theta_0)) = b(t_n;\theta=\theta_0) + dt \cdot \sum_{i=n}^{N-1} h_{t_i},$$

then

$$V^C(b(t_n;\theta=\theta_0)+C,t_n;\theta=\theta_0+C,b(t_n)=b(t_n;\theta=\theta_0)+C) = b(t_n;\theta=\theta_0) + C + dt \cdot \sum_{i=n}^{N-1} h_{t_i}.$$

Thus, we get that

$$b(t_n;\theta=\theta_0+C) = b(t_n;\theta=\theta_0) + C, \text{ for } 0 \leq t \leq T.$$

∎

Figure 1 exemplifies the threshold function, $b(t)$ for a time horizon of $T=15$, with different $\theta$ values.

**Lemma 5** $b(t_n;\mathbf{h}=\mathbf{h}^H) \leq b(t_n;\mathbf{h}=\mathbf{h}^L)$ where $\mathbf{h}^H = \{h_{t_0}^H,\ldots,h_{t_{N-1}}^H\}$, $\mathbf{h}^L = \{h_{t_0}^L,\ldots,h_{t_{N-1}}^L\}$, $h_{t_i}^L \leq h_{t_i}^H$, for $n \leq i \leq N-1$ and $1 \leq n < N-1$.

**Proof.** Let $\mathbf{b}^L = \{b^L(t_0),\ldots,b^L(t_N)\}$ denote the optimal set of threshold values for a system with holding cost $\mathbf{h}^L$, and let $V^C(x,t_n;\mathbf{b}=\mathbf{b}^L,\mathbf{h}=\mathbf{h}^L)$ denote the value function from *continuation* decision for a system with holding cost $\mathbf{h}^L$, when following thresholds $\mathbf{b}^L$ for $n+1 \leq i \leq N-1$. According to eq. (6),

$$b^L(t_n) + dt \cdot \sum_{i=n}^{N-1} h_{t_i}^L = V^C\left(b^L(t_n),t_n;\mathbf{b}=\mathbf{b}^L,\mathbf{h}=\mathbf{h}^L\right). \quad (32)$$

Let $h_{t_i}^H \geq h_{t_i}^L$ for $n \leq i \leq N-1$, and let $V^C(x,t_n;\mathbf{b}=\mathbf{b}^L,\mathbf{h}=\mathbf{h}^H)$ denote the value function from *continuation* decision for a system with holding cost $\mathbf{h}^H$, when following thresholds $b^L(t_i)$ for $n+1 \leq i \leq N-1$. Under this setting, according to eq. (8) we get that the difference between $V^C\left(b^L(t_n),t_n;\mathbf{b}=\mathbf{b}^L,\mathbf{h}=\mathbf{h}^L\right)$

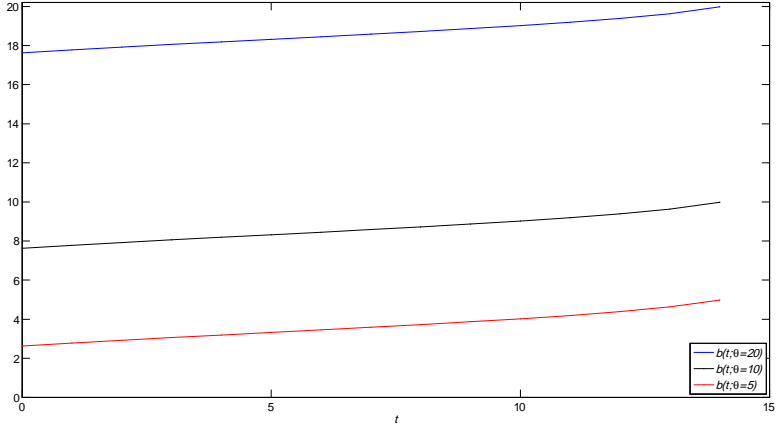

Figure 1: Threshold function, $b(t, \theta = \theta_0)$ for $T = 15$, $dt = 1$, $K = 0.5$, $\sigma = 1$, $h_t = (T - t) \cdot 10^{-2}$, for $\theta = 20, 10, 5$.

and $V^C\left(b^L(t_n), t_n; \mathbf{b} = \mathbf{b}^L, \mathbf{h} = \mathbf{h}^H\right)$ is only in the holding cost. That is,

$$V^C\left(b^L(t_n), t_n; \mathbf{b} = \mathbf{b}^L, \mathbf{h} = \mathbf{h}^L\right) = V^C\left(b^L(t_n), t_n; \mathbf{b} = \mathbf{b}^L, \mathbf{h} = \mathbf{h}^H\right)$$
$$- \sum_{i=n+1}^{N-1}\left[P_{t_n, t_i}(b^L(t_n); \mathbf{b} = \mathbf{b}^L) \cdot dt \cdot \sum_{j=i}^{N-1}\left(h_{t_j}^H - h_{t_j}^L\right)\right]. \tag{33}$$

Then, from eqs. $(32), (33)$ we get that,

$$b^L(t_n) + dt \cdot \sum_{i=n}^{N-1} h_{t_i}^L + \sum_{i=n+1}^{N-1}\left[P_{t_n, t_i}(b^L(t_n); \mathbf{b} = \mathbf{b}^L) \cdot dt \cdot \sum_{j=i}^{N-1}\left(h_{t_j}^H - h_{t_j}^L\right)\right]$$
$$= V^C\left(b^L(t_n), t_n; \mathbf{b} = \mathbf{b}^L, \mathbf{h} = \mathbf{h}^H\right). \tag{34}$$

Moreover,

$$b^L(t_n) + dt \cdot \sum_{i=n}^{N-1} h_{t_i}^L + \sum_{i=n+1}^{N-1}\left[P_{t_n, t_i}(b^L(t_n); \mathbf{b} = \mathbf{b}^L) \cdot dt \cdot \sum_{j=i}^{N-1}\left(h_{t_j}^H - h_{t_j}^L\right)\right] \tag{35}$$
$$\leq b^L(t_n) + dt \cdot \sum_{i=n}^{N-1} h_{t_i}^H,$$

and
$$V^C\left(b^L(t_n), t_n; \mathbf{b} = \mathbf{b}^H, \mathbf{h} = \mathbf{h}^H\right) \leq V^C\left(b^L(t_n), t_n; \mathbf{b} = \mathbf{b}^L, \mathbf{h} = \mathbf{h}^H\right), \tag{36}$$

as $V^C\left(b^L(t_n), t_n; \mathbf{b} = \mathbf{b}^H, \mathbf{h} = \mathbf{h}^H\right)$ represents the optimal value from *continuation* decision for a system with holding cost $h_{t_i}^H$. Therefore, we get from eqs. $(34), (35), (36)$ that

$$V^C\left(b^L(t_n), t_n; \mathbf{b} = \mathbf{b}^H, \mathbf{h} = \mathbf{h}^H\right) \leq b^L(t_n) + dt \cdot \sum_{i=n}^{N-1} h_{t_i}^H.$$

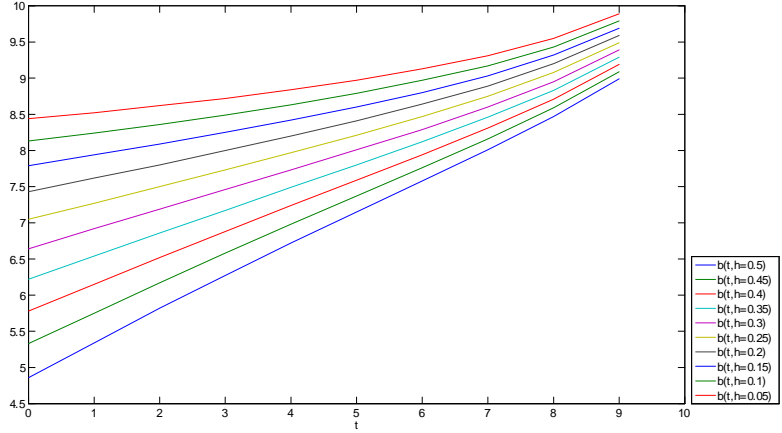

Figure 2: Threshold function, $b(t)$ for $T = 10$, $dt = 1$, $K = 0.5$, $\sigma = 1$, $\theta = 10$ across varying holding values, $h$.

Hence, for a system with holding cost of $\mathbf{h}^H$, it is better to wait with the purchasing price at price $b^L(t_n)$, which together with the uniqueness of the threshold value (Lemma 1) implies that $b\left(t_n; \mathbf{h} = \mathbf{h}^H\right) \leq b\left(t_n; \mathbf{h} = \mathbf{h}^L\right)$ for $h_{t_i}^L \leq h_{t_i}^H, n \leq i \leq N-1$ and $1 \leq n < N-1$. ∎

Figure 2 exemplifies the threshold function, $b(t)$ for a time horizon of $T = 10$, with different holding values.

## 4.1 Threshold function properties under special cases of $h$

In this section, we analyze additional properties of the threshold function when the holding cost is subject to follow specific forms. We denote the holding cost as a function of a set of parameters as $h(\cdot)$. In addition, for better readability of the proofs we assume in this section that $\theta = 0$. Note that according to Lemma 4 the results can be easily adjusted to a general $\theta$, as $\theta$ only affects the threshold by a constant.

**Lemma 6** $b(t_n)$ is linear to $\sigma$, under an holding cost that is linear to $\sigma$, $h_{t_i}(\sigma)$, for $n \leq i \leq N-1$.

**Proof.** We prove this Lemma by induction. First, we note that under the holding cost, $h_{t_{N-1}}(\sigma)$, we get that $b(t_{N-1})$ is linear to $\sigma$: according to Corollary 1, $b(t_{N-1}) = -h_{t_{N-1}}(\sigma)/K$ (under $\theta = 0$) which is linear to $\sigma$. Under the induction assumption, $b(t_i)$ is linear to $\sigma$, under an holding cost, $h_{t_i}(\sigma)$, for $t_{n+1} \leq t_i \leq t_{N-1}$ (Assumption 1).

Next we show that under Assumption 1 with $h_{t_i}(\sigma)$, $b(t_n)$ is also linear to $\sigma$. We show this by proving that under Assumption 1 the following Properties hold:

1. According to eq. (14), we get that $\text{cov}\left(Z_{t_n,t_l}(x), Z_{t_n,t_j}(x)\right)$ is independent of $\sigma$ and $h_{t_i}(\sigma)$, for $n+1 \leq l \leq N, n < j \leq N$.

2. According to eq. (15), for $\sigma = \sigma_0$ and $X_{t_n} = x$, we get that

$$\beta_{t_n,t_l}(x; \sigma = \sigma_0) = \frac{b(t_l; \sigma = \sigma_0) - \left[x \cdot (1-dtK)^{l-n}\right]}{\sigma_{t_n,t_l}(;\sigma = \sigma_0)}.$$

For a constant $C$, we get from Assumption 1 and eq. (11) that

$$b(t_l; \sigma = C \cdot \sigma_0)/\sigma_{t_n,t_l}(;\sigma = C \cdot \sigma_0) = b(t_l; \sigma = \sigma_0)/\sigma_{t_n,t_l}(;\sigma = \sigma_0) \text{ for } n+1 \leq l \leq N. \quad (37)$$

Then, for $X_{t_n} = C \cdot x$, we get

$$\frac{b\left(t_l; \sigma = C \cdot \sigma_0\right) - \left[C \cdot x \cdot (1 - dtK)^{l-n}\right]}{\sigma_{t_n, t_l}\left(; \sigma = C \cdot \sigma_0\right)} = \frac{b\left(t_l; \sigma = \sigma_0\right) - \left[x \cdot (1 - dtK)^{l-n}\right]}{\sigma_{t_n, t_l}\left(; \sigma = \sigma_0\right)} = \beta_{t_n, t_l}\left(x; \sigma = \sigma_0\right).$$

Hence,

$$\beta_{t_n, t_l}\left(x; \sigma = \sigma_0\right) = \beta_{t_n, t_l}\left(C \cdot x; \sigma = C \cdot \sigma_0\right).$$

3. According to Properties 1 and 2 we get that

$$F\left(-\hat{\mathbf{B}}_{t_n, t_i}(x), \Sigma_{t_n, t_i}; \sigma = \sigma_0\right) = F\left(-\hat{\mathbf{B}}_{t_n, t_i}(C \cdot x), \Sigma_{t_n, t_i}; \sigma = C \cdot \sigma_0\right),$$

and

$$F\left(-\mathbf{B}_{t_n, t_i}(x), \Sigma_{t_n, t_i}; \sigma = \sigma_0\right) = F\left(-\mathbf{B}_{t_n, t_i}(C \cdot x), \Sigma_{t_n, t_i}; \sigma = C \cdot \sigma_0\right),$$
$$\text{for } t_{n+1} \leq t_i \leq t_N.$$

Then, according to eq. (16),

$$P_{t_n, t_i}\left(x; \sigma = \sigma_0\right) = P_{t_n, t_i}\left(C \cdot x; \sigma = C \cdot \sigma_0\right).$$

4. According to Properties $1 - 3$ and eqs. $(17) - (23)$, we get that

$$E\left(Z_{t_n, t_i}(x) | \tau_{t_n}(x) > t_i; \sigma = \sigma_0\right) = E\left(Z_{t_n, t_i}(C \cdot x) | \tau_{t_n}(C \cdot x) > t_i; \sigma = C \cdot \sigma_0\right),$$

and

$$E\left(Z_{t_n, t_i}(x) | \tau_{t_n}(x) > t_{i-1}; \sigma = \sigma_0\right) = E\left(Z_{t_n, t_i}(C \cdot x) | \tau_{t_n}(C \cdot x) > t_{i-1}; \sigma = C \cdot \sigma_0\right).$$

Then, according to eqs. (9), (11), (24), and (25) we get that

$$C \cdot E\left(X_{t_n, t_i}(x) | \tau_{t_n}(x) > t_i; \sigma = \sigma_0\right) = E\left(X_{t_n, t_i}(C \cdot x) | \tau_{t_n}(C \cdot x) > t_i; \sigma = C \cdot \sigma_0\right),$$

and

$$C \cdot E\left(X_{t_n, t_i}(x) | \tau_{t_n}(x) > t_{i-1}; \sigma = \sigma_0\right) = E\left(X_{t_n, t_i}(C \cdot x) | \tau_{t_n}(C \cdot x) > t_{i-1}; \sigma = C \cdot \sigma_0\right),$$

which implies according to eq. (27) that

$$C \cdot E_{t_n, t_i}\left(x; \sigma = \sigma_0\right) = E_{t_n, t_i}\left(C \cdot x; \sigma = C \cdot \sigma_0\right),$$

for $t_{n+1} \leq t_i \leq t_N$.

5. According to Properties 3, 4 and eq. (8), we get that in the continuation region

$$C \cdot V^C\left(x, t_n; \sigma = \sigma_0\right) = C \cdot \sum_{i=n+1}^{N} P_{t_n, t_i}\left(x; \sigma = \sigma_0\right) \cdot \left(E_{t_n, t_i}\left(x; \sigma = \sigma_0\right) + dt \cdot \sum_{j=i}^{N-1} h_{t_j}(\sigma_0)\right)$$

$$= \sum_{i=n+1}^{N} P_{t_n, t_i}\left(C \cdot x; \sigma = C \cdot \sigma_0\right) \cdot \left(E_{t_n, t_i}\left(C \cdot x; \sigma = C \cdot \sigma_0\right) + dt \cdot \sum_{j=i}^{N-1} h_{t_j}(C \cdot \sigma_0)\right),$$

and therefore

$$C \cdot V^C\left(x, t_n; \sigma = \sigma_0\right) = V^C\left(C \cdot x, t_n; \sigma = C \cdot \sigma_0\right), \tag{38}$$

for $h_{t_i}(\sigma)$ that is linear to $\sigma$, for $t_{n+1} \leq t_i \leq t_N$.

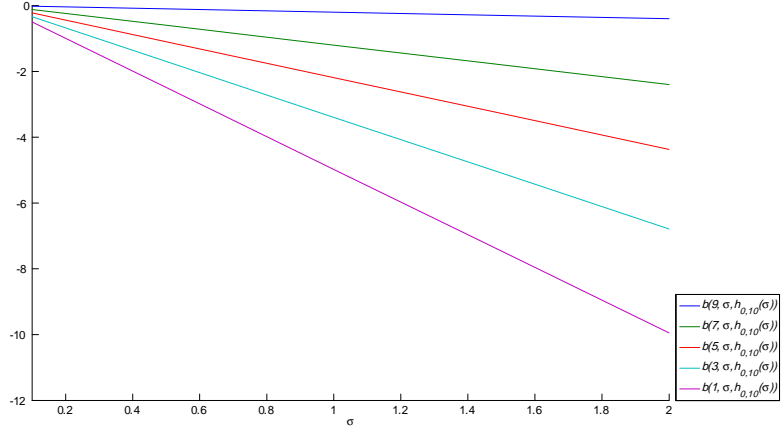

Figure 3: Threshold function, $b(t_n; \sigma, h_{0,t_{10}})$ for $t_n = 1, 3, 5, 7, 9, T = 10, dt = 1, K = 0.5$, for $0.1 \leq \sigma \leq 2$, and $h_{t_i} = \sigma \cdot (t_N - t_i) \cdot 10^{-1}$.

Finally, according to eqs. $(6), (38)$, under Assumption 1, for a given threshold, $b(t_n; \sigma = \sigma_0)$

$$C \cdot \left( b(t_n; \sigma = \sigma_0) + dt \cdot \sum_{i=n}^{N-1} h_{t_i}(\sigma_0) \right) = V^C \left( C \cdot b(t_n; \sigma = \sigma_0), t_n; \sigma = C \cdot \sigma_0 \right)$$

which implies that under Assumption 1,

$$C \cdot b(t_n; \sigma = \sigma_0) = b(t_n; \sigma = C \cdot \sigma_0).$$

Hence, $b(t_n)$ is linear to $\sigma$, under $h_{t_i}(\sigma)$, for $0 \leq t_n \leq t_{N-1}$. ∎

Figure 3 exemplifies the threshold function, $b(t; \sigma)$ for $t = 1, 3, 5, 7, 9, T = 10, dt = 1, K = 0.5$, and $h_{t_i}(\sigma) = \sigma \cdot (t_N - t_i) \cdot 10^{-1}$.

**Lemma 7** $b(t_n)$ is linear to $\sqrt{dt}$ and to $1/\sqrt{K}$, under a constant product, $dt \cdot K$, and a holding cost that is linear to $\sqrt{K}$ and $1/\sqrt{dt}$, $h_{t_i}(dt, K)$, for $n \leq i \leq N-1$.

**Proof.** Similar to the proof in Lemma 6, it can be proved that $b(t_n)$ is linear to $\sqrt{dt}$ and to $1/\sqrt{K}$, under a constant product of $dt \cdot K$, and holding cost, $h_{t_i}(dt, K)$. In order to adjust the proof in Lemma 6, one should follow the same steps as in the proof of Lemma 6, with the following adjustments:

- $h_{t_i}(dt, K)$ replacing $h_{t_i}(\sigma)$.
- $\sqrt{dt}$ replacing $\sigma$.
- $dt = dt_0, K = K_0$ replacing $\sigma = \sigma_0$.
- $dt = C \cdot dt_0, K = K_0/C$ replacing $\sigma = C \cdot \sigma_0$.
- $\sqrt{C} \cdot x$ replacing $C \cdot x$.
- $\sqrt{C} \cdot E\left(X_{t_n,t_i}(x) | \tau_{t_n}(x) > t_i; dt = dt_0, K = K_0\right)$ replacing $C \cdot E\left(X_{t_n,t_i}(x) | \tau_{t_n}(x) > t_i; \sigma = \sigma_0\right)$
- $\sqrt{C} \cdot E\left(X_{t_n,t_i}(x) | \tau_{t_n}(x) > t_{i-1}; dt = dt_0, K = K_0\right)$ replacing $C \cdot E\left(X_{t_n,t_i}(x) | \tau_{t_n}(x) > t_{i-1}; \sigma = \sigma_0\right)$
- $\sqrt{C} \cdot E_{t_n,t_i}(x; dt = dt_0, K = K_0)$ replacing $C \cdot E_{t_n,t_i}(x; \sigma = \sigma_0)$.

- $$\sqrt{C} \cdot V^C(x, t_n; dt = dt_0, K = K_0)$$
$$= \sqrt{C} \cdot \sum_{i=n+1}^{N} P_{t_n, t_i}(x; dt = dt_0, K = K_0) \cdot \left( E_{t_n, t_i}(x; dt = dt_0, K = K_0) + dt_0 \cdot \sum_{j=i}^{N-1} h_{t_j}(dt_0, K_0) \right)$$
$$= \sum_{i=n+1}^{N} P_{t_n, t_i}\left(\sqrt{C} \cdot x; dt = C \cdot dt_0, K = K_0/C\right) \cdot \left( \begin{array}{c} E_{t_n, t_i}\left(\sqrt{C} \cdot x; dt = C \cdot dt_0, K = K_0/C\right) \\ + C \cdot dt_0 \cdot \sum_{j=i}^{N-1} h_{t_j}(C \cdot dt_0, K_0/C) \end{array} \right),$$

  replacing
$$C \cdot V^C(x, t_n; \sigma = \sigma_0) = C \cdot \sum_{i=n+1}^{N} P_{t_n, t_i}(x; \sigma = \sigma_0) \cdot \left( E_{t_n, t_i}(x; \sigma = \sigma_0) + dt \cdot \sum_{j=i}^{N-1} h_{t_j}(\sigma_0) \right)$$
$$= \sum_{i=n+1}^{N} P_{t_n, t_i}(C \cdot x; \sigma = C \cdot \sigma_0) \cdot \left( E_{t_n, t_i}(C \cdot x; \sigma = C \cdot \sigma_0) + dt \cdot \sum_{j=i}^{N-1} h_{t_j}(C \cdot \sigma_0) \right),$$

- $$\sqrt{C} \cdot V^C(x, t_n; dt = dt_0, K = K_0) = V^C\left(\sqrt{C} \cdot x, t_n; dt = C \cdot dt_0, K = K_0/C\right)$$

  replacing
$$C \cdot V^C(x, t_n; \sigma = \sigma_0) = V^C(C \cdot x, t_n; \sigma = C \cdot \sigma_0)$$

- $$\sqrt{C} \cdot \left( b(t_n; dt = dt_0, K = K_0) + dt_0 \cdot \sum_{i=n}^{N-1} h_{t_i}(dt_0, K_0) \right)$$
$$= V^C\left(\sqrt{C} \cdot b(t_n; dt = dt_0, K = K_0), t_n; dt = C \cdot dt_0, K = K_0/C\right)$$

  replacing
$$C \cdot \left( b(t_n; \sigma = \sigma_0) + dt \cdot \sum_{i=n}^{N-1} h_{t_i}(\sigma_0) \right) = V^C(C \cdot b(t_n; \sigma = \sigma_0), t_n; \sigma = C \cdot \sigma_0)$$

- $$\sqrt{C} \cdot b(t_n; dt = dt_0, K = K_0) = b(t_n; dt = C \cdot dt_0, K = K_0/C)$$

  replacing
$$C \cdot b(t_n; \sigma = \sigma_0) = b(t_n; \sigma = C \cdot \sigma_0).$$

Under this adjustments, the same induction process as in the proof of Lemma 6 derives that $b(t_n)$ is linear to $\sqrt{dt}$, under a constant product, $dt \cdot K$, and an holding cost that is linear to $\sqrt{K}$ and $1/\sqrt{dt}$. Note that as $dt \cdot K$ is constant, than we get that $b(t_n)$ is linear to $1/\sqrt{K}$ as well. ∎

Figure 4 exemplifies the threshold function, $b(t_n; dt, K)$ for $n = 1, 3, 5, 7, 9, N = 10$, $dt \cdot K = 0.8$ and $h_{t_i}(dt, K) = 1/\sqrt{dt} \cdot (t_N - t_i) \cdot 10^{-1}$.

## 5 Conclusion

In this research, we develop an optimal policy for managing inventory of an item with a stochastic mean reverting price process. Our inventory problem considers a discrete finite time horizon in which the item has to be bought, and a convex holding cost. The item's price fluctuates according to an $AR(1)$ process, that reflects a mean-reversion behavior. Under this setting we prove that the optimal policy is in the form of a unique threshold function, that determines for each time point, the price region for which a purchasing decision is optimal. We provide a simple algorithm to calculate the threshold function. This algorithm is based on the Bellman equation that is used to calculate the value function of the system in each step. The value function calculation relies on the *Crossing Time Probability* and *Overshoot Expectation*, that are given in a closed form solution. We analyze the threshold function with respect to the model's logistic parameters, and to the price process parameters. The analysis derives interesting properties of the threshold function, that are proved theoretically and exemplified using graphs.

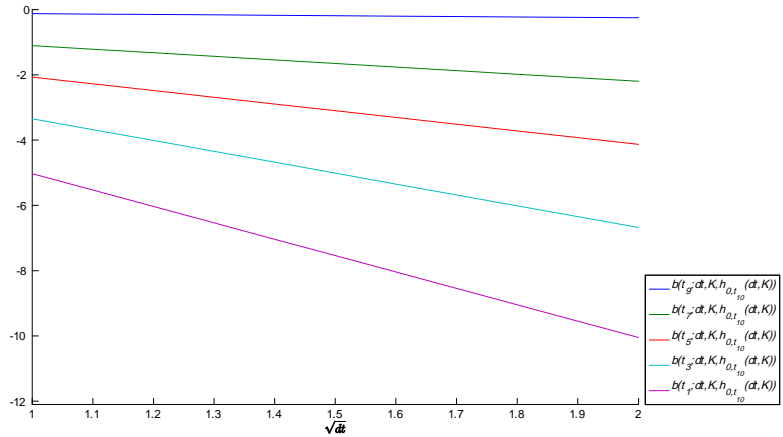

Figure 4: Threshold function, $b\left(t_n; dt, K, h_{0,t_{10}}\right)$ for $n = 1, 3, 5, 7, 9, N = 10$, $dt \cdot K = 0.8$, $\sigma = 1$, for $1 \leq \sqrt{dt} \leq 2$, and $h_{t_i} = 1/\sqrt{dt} \cdot (t_N - t_i) \cdot 10^{-1}$.

# References

[1] Alili, L., Patie, P., and Pedersen, J. L. (2005). Representations of the first hitting time density of an Ornstein-Uhlenbeck process 1. *Stochastic Models*, **21(4)**, 967-980.

[2] Bellman, R. E., and Dreyfus, S. E. (2015). Applied dynamic programming. Princeton university press. Chicago.

[3] Berling, P. and Martinez-de-Albeniz, V. (2011). Optimal inventory policies when purchase price and demand are stochastic. *Operations Research*, **59(1)**, 109-124.

[4] Bernard, J.T., Khalaf, L., Kichian, M. and McMahon, S. (2008). Forecasting commodity prices: GARCH, jumps, and mean reversion. *Journal of Forecasting*, **27**, 279-291.

[5] Bessembinder, H., Coughenour, J.F, Sequin, P.J, and Smoller, M.M (1995). Mean reversion in equilibrium asset prices: evidence from futures term structure. *The Journal of Finance*, **50(1)**, 361-375.

[6] Black, F., and Scholes, M. (1973). The pricing of options and corporate liabilities. *Journal of Political Economy,* **81(3)**, 637-659.

[7] Christensen, S. (2012). Phase-type distributions and optimal stopping for autoregressive processes. *Journal of Applied Probability*, 22-39.

[8] Christensen, S., Irle, A., and Novikov, A. (2011). An elementary approach to optimal stopping problems for AR (1) sequences. *Sequential Analysis*, **30(1)**, 79-93.

[9] Elliott, R. J., Van Der Hoek, J., and Malcolm, W. P. (2005). Pairs trading. *Quantitative Finance*, 5(3), 271-276.

[10] Finster, M. (1982). Optimal stopping on autoregressive schemes. *The Annals of Probability*, **10(3)**, 745-753.

[11] Gavirneni, S. (2004). Periodic review inventory control with fluctuating purchasing costs. *Operations Research Letters*, **32(4)**, 374-379.

[12] Gil-Alana, L. A. (2000). Mean reversion in the real exchange rates. *Economics Letters*, **69(3)**, 285-288.

[13] Goel, A., and Gutierrez, G. J. (2011). Multiechelon procurement and distribution policies for traded commodities. *Management Science*, **57(12)**, 2228-2244.

[14] Golabi, K. (1985). Optimal inventory policies when ordering prices are random. *Operations Research*, **33(3)**, 575-588.

[15] Hahn, W. J., and Dyer, J. S. (2011). A discrete time approach for modeling two-factor mean-reverting stochastic processes. *Decision Analysis*, **8(3)**, 220-232.

[16] Haksöz, Ç. and Seshadri, S. (2007). Supply chain operations in the presence of a spot market: A review with discussion. *J. Oper. Res. Soc.* **58(11)**, 1412-1429.

[17] Higgs, H., and Worthington, A. (2008). Stochastic price modeling of high volatility, mean-reverting, spike-prone commodities: The Australian wholesale spot electricity market. *Energy Economics*, **30(6)**, 3172-3185.

[18] Hull, J.C. (2011). *Options, Futures and Other Derivatives*, 8th ed. Pearson Education.

[19] Hull, J. C., and White, A. D. (1994). Numerical procedures for implementing term structure models I: Single-factor models. *The Journal of Derivatives*, **2(1)**, 7-16.

[20] Jacka, S. 1. (1991). Optimal stopping and the American put. *Mathematical Finance*, **1(2)**, 1-14.

[21] Kalymon, B.A. (1971). Stochastic prices in a single-item inventory purchasing model. *Operations Research*, **19(6)**, 1434-1458.

[22] Li, C. L., & Kouvelis, P. (1999). Flexible and risk-sharing supply contracts under price uncertainty. *Management Science*, **45(10)**, 1378-1398.

[23] Longstaff, F. A., and Schwartz, E. S. (2001). Valuing American options by simulation: a simple least-squares approach. *Review of Financial studies*, **14(1)**, 113-147.

[24] Nelson, D. B., and Ramaswamy, K. (1990). Simple binomial processes as diffusion approximations in financial models. *Review of Financial Studies*, **3(3)**, 393-430.

[25] Novikov, A. A. (2009). On distributions of first passage times and optimal stopping of AR (1) sequences. *Theory of Probability & Its Applications*, **53(3)**, 419-429.

[26] Pindyck, R.S. (2001). The dynamics of commodity spot and futures. *The Energy Journal*, **22(3)**, 1-29.

[27] Pindyck, R.S. and Rubinfeld, D.L. (1998). *Econometric Models and Economic Forecasts*, 4th edn McGraw-Hill, 1998.

[28] Schwartz, E.S (1997). The stochastic behavior of commodity prices: implications for valuation and hedging. *The Journal of Finance,* **5(3)**, 923-973.

[29] Tallis, G. M. (1961). The moment generating function of the truncated multi-normal distribution, *Journal of the Royal Statistical Society. Series B (Methodological)*, **23(1)**, 23-29.

[30] Uhlenbeck, G. E., and Ornstein, L. S. (1930). On the theory of the Brownian motion. *Physical review*, **36(5)**, 823.

[31] Wang, Y. (2001). The optimality of myopic stocking policies for systems with decreasing purchasing prices. *European Journal of Operational Research*, **133(1)**, 153-159.

[32] William, H. G. (2003). Econometric Analysis, Prentice Hall, Upper Saddle River, NJ.

[33] Zhang, K., Nieto, A., and Kleit, A.N (2015). The real option value of mining operations using mean-reverting commodity prices. *Mineral Economics*, **28**, 11-22.

**Appendix 1** *Proof of Lemma 1.*

**Proof.** First we prove Property 1: at any time point, $0 \leq t_n \leq t_{N-1}$, there exists a finite price, $x^H(t_n) \equiv \theta - h_{t_{N-1}}/K$, for which any $x \geq x^H(t_n)$ satisfies $x + dt \cdot \sum_{i=n}^{N-1} h_{t_i} > V^C(x, t_n)$. A feasible policy at any finite price $x > \theta - h_{t_n}/K$ is to wait with the purchasing decision at time $t_n$, and purchase the item immediately at the next time interval at price $X_{t_n, t_{n+1}}(x)$. According to eq. (9), this decision leads to expected cost of:

$$E\left[X_{t_n, t_{n+1}}(x)\right] + dt \cdot \sum_{i=n+1}^{N-1} h_{t_i}$$
$$= \theta + (x - \theta)(1 - dt \cdot K) + dt \cdot \sum_{i=n+1}^{N-1} h_{t_i}, \tag{39}$$

where for any $x > \theta - h_{t_n}/K$, this expectation is lower than the cost of of purchasing the item at time $t_n$. That is,

$$x + dt \cdot \sum_{i=n}^{N-1} h_{t_i} > \theta + (x - \theta)(1 - dt \cdot K) + dt \cdot \sum_{i=n+1}^{N-1} h_{t_i}, \text{ for any } x > \theta - h_{t_n}/K, \tag{40}$$

where the *left hand side (LHS)* represents the cost of a purchasing decision at time $t_n$, and the *right hand side (RHS)* represents the expected cost under a purchasing decision at time $t_{n+1}$. By definition, the optimal policy has a lower cost than any feasible policy, and therefore we get that $E\left[X_{t_n, t_{n+1}}(x)\right] + dt \cdot \sum_{i=n+1}^{N-1} h_{t_i} \geq V^C(x, t_n)$. Hence, according to eqs. (40), and (39), $x + dt \cdot \sum_{i=n}^{N-1} h_{t_i} > V^C(x, t_n)$, for $0 \leq t_n < t_{N-1}$, for any $x > \theta - h_{t_n}/K$.

Next we prove by induction that Properties 2, 3 hold, and there exists a finite single threshold value, $b(t_n)$, for $0 \leq t_n < t_{N-1}$:

Property 2 - At any time point, $0 \leq t_n \leq t_{N-1}$, there exists a finite price, $x^L(t_n) \equiv \min\left\{\frac{b(t_{n+1}) - \theta dt \cdot K}{(1 - dt \cdot K)}, \frac{\theta \cdot dt \cdot K - \sigma \cdot \sqrt{dt} - h_{t_n} \cdot dt}{dt \cdot K}\right\}$, for which any $x \leq x^L(t_n)$ satisfies $x + dt \cdot \sum_{i=n}^{N} h_i \leq V^C(x, t_n)$.

Property 3 - $V^C(x, t_n)$ is concave and increasing in $x$.

First we verify that for time $t_{N-1}$ a single threshold exists, and that Properties 2 and 3 hold: According to Corollary 1, at time $t_{N-1}$, a single threshold exists at the value of $b(t_{N-1}) = \theta - h_{t_{N-1}}/K$, and any price $x < \theta - h_{t_{N-1}}/K$ satisfies $x + dt \cdot \sum_{i=n}^{N} h_i < V^C(x, t_n)$. Note that $\theta - h_{t_{N-1}}/K > \frac{\theta \cdot dt \cdot K - \sigma \cdot \sqrt{dt} - h_{t_{N-1}} \cdot dt}{dt \cdot K}$, and therefore Property 2 holds. In addition, at time $t_{N-1}$ the value function under continuation decision leads to a definite purchase at the end of the planning horizon at time $t_N$. According to eq. (9)

$$V^C(x, t_{N-1}) = E\left[X_{t_{N-1}, t_N}(x)\right] = \theta + (x - \theta)(1 - dt \cdot K).$$

That is, $V^C(x, t_{N-1})$ is concave (linear) and increasing in $x$. Hence, Properties 3 also holds for time $t_{N-1}$.

Next we assume that for time $t_j : t_{n+1} \leq t_j \leq t_N$, there exists a finite single threshold value, $b(t_j)$, and that Properties 2 and 3 hold. Under this assumption, we prove by induction that the above also holds for time $t_n$. From Property 1, we get that

$$V^C(x, t_n) = \Pr\left(X_{t_n, t_{n+1}}(x) > b(t_{n+1})\right) \cdot E\left[V^C(X_{t_n, t_{n+1}}(x), t_{n+1}) | X_{t_n, t_{n+1}}(x) > b(t_{n+1})\right] \quad (41)$$
$$+ \left(1 - \Pr\left(X_{t_n, t_{n+1}}(x) > b(t_{n+1})\right)\right) \cdot \left(E\left[X_{t_n, t_{n+1}}(x) | X_{t_n, t_{n+1}}(x) \leq b(t_{n+1})\right] + dt \cdot \sum_{i=n+1}^{N-1} h_{t_i}\right).$$

Note that according to threshold policy,

$$E\left[V^C(X_{t_n, t_{n+1}}(x), t_{n+1}) | X_{t_n, t_{n+1}}(x) > b(t_{n+1})\right] \geq E\left[X_{t_n, t_{n+1}}(x) | X_{t_n, t_{n+1}}(x) \leq b(t_{n+1})\right] + dt \cdot \sum_{i=n+1}^{N-1} h_{t_i}.$$

Therefore, we can substitute the *RHS* of eq. (41), and get the following inequality:

$$V^C(x, t_n) \geq E\left[X_{t_n, t_{n+1}}(x) | X_{t_n, t_{n+1}}(x) \leq b(t_{n+1})\right] + dt \cdot \sum_{i=n+1}^{N-1} h_{t_i}. \quad (42)$$

According to the induction assumption, and eq. (9) there exists a finite $\bar{x} \equiv \frac{b(t_{n+1}) - \theta dt \cdot K}{(1 - dt \cdot K)}$ that satisfies $E\left[X_{t_n, t_{n+1}}(\bar{x})\right] = b(t_{n+1})$. Hence, we can infer from eq. (9) that $E\left[X_{t_n, t_{n+1}}(x)\right] \leq b(t_{n+1}) \iff x \leq \bar{x}$. Therefore, we get that for $x \leq \bar{x}$:

$$E\left[X_{t_n, t_{n+1}}(x) | X_{t_n, t_{n+1}}(x) \leq b(t_{n+1})\right] + dt \cdot \sum_{i=n+1}^{N-1} h_{t_i}$$
$$\geq E\left[X_{t_n, t_{n+1}}(x) | X_{t_n, t_{n+1}}(x) \leq E\left[X_{t_n, t_{n+1}}(x)\right]\right] + dt \cdot \sum_{i=n+1}^{N-1} h_{t_i}, \quad (43)$$

where by eq. (3)

$$E\left[X_{t_n, t_{n+1}}(x) | X_{t_n, t_{n+1}}(x) \leq E\left[X_{t_n, t_{n+1}}(x)\right]\right] + dt \cdot \sum_{i=n+1}^{N-1} h_{t_i}$$
$$= E\left[\theta + (x - \theta)(1 - dt \cdot K) + \sigma \cdot \varepsilon_{t_{n+1}} | \varepsilon_{t_{n+1}} \leq 0\right] + dt \cdot \sum_{i=n+1}^{N-1} h_{t_i}$$
$$= \theta + (x - \theta)(1 - dt \cdot K) + \sigma \cdot E\left[\varepsilon_{t_{n+1}} | \varepsilon_{t_{n+1}} \leq 0\right] + dt \cdot \sum_{i=n+1}^{N-1} h_{t_i}. \quad (44)$$

The term, $\varepsilon_{t_{n+1}} | \varepsilon_{t_{n+1}} < 0$, distributes as a normal variable, with zero mean and $dt$ variance, that is truncated below its mean. This term has a known expected value of (see: [32]):

$$E\left[\varepsilon_{t_{n+1}} | \varepsilon_{t_{n+1}} \leq 0\right] = -\sqrt{dt}\phi(0)/\Phi(0), \quad (45)$$

where $\phi(\cdot)$ and $\Phi(\cdot)$ denotes the *PDF* and the cumulative distribution function (*CDF*) respectively, of a standard normal random variable. As $0 < \phi(0)/\Phi(0) < 1$, and $\sqrt{dt}, \sigma > 0$, we get from eqs. (42,43,44,45) that

$$V^C(x, t_n) > \theta + (x - \theta)(1 - dt \cdot K) - \sigma \cdot \sqrt{dt} + dt \cdot \sum_{i=n+1}^{N-1} h_{t_i}, \text{ for } x \leq \frac{b(t_{n+1}) - \theta dt \cdot K}{(1 - dt \cdot K)}, \quad (46)$$

where,

$$\theta + (x - \theta)(1 - dt \cdot K) - \sigma \cdot \sqrt{dt} + dt \cdot \sum_{i=n+1}^{N-1} h_{t_i} > x + dt \cdot \sum_{i=n}^{N-1} h_{t_i}$$
$$\iff x < \frac{\theta \cdot dt \cdot K - \sigma \cdot \sqrt{dt} - h_{t_n} \cdot dt}{dt \cdot K}. \tag{47}$$

Hence, by eqs. (46,47) we get that there exists a finite price, $x^L(t_n) \equiv \min\left\{\frac{b(t_{n+1}) - \theta dt \cdot K}{(1 - dt \cdot K)}, \frac{\theta \cdot dt \cdot K - \sigma \cdot \sqrt{dt} - h_{t_n} \cdot dt}{dt \cdot K}\right\}$, for which any $x \leq x^L(t_n)$ satisfies:

$$V^C(x, t_n) > x + dt \cdot \sum_{i=n}^{N-1} h_{t_i}.$$

That is, Property 2 is valid for time $t_n$. Next, we note that under the induction assumption there exists a finite single threshold value, $b(t_{n+1})$, for time $t_{n+1}$ and therefore we get that

$$V^C(x, t_n) = E\left[\min\left\{X_{t_n, t_{n+1}}(x) + dt \cdot \sum_{i=n+1}^{N-1} h_{t_i}, V^C(X_{t_n, t_{n+1}}(x), t_{n+1})\right\}\right]. \tag{48}$$

We define $\varepsilon$, as a realization of $\varepsilon_{t_{n+1}}$ under its probability density function, $g(\varepsilon)$, and $x_{t_n, t_{n+1}}(x, \varepsilon)$, as the value of $X_{t_n, t_{n+1}}(x)$ when $\varepsilon_{t_{n+1}} = \varepsilon$. Therefore, we can set eq. (48) as

$$V^C(x, t_n) = \int_{-\infty}^{\infty} \min\left\{x_{t_n, t_{n+1}}(x, \varepsilon) + dt \cdot \sum_{i=n+1}^{N-1} h_{t_i}, V^C(x_{t_n, t_{n+1}}(x, \varepsilon), t_{n+1})\right\} g(\varepsilon) d\varepsilon. \tag{49}$$

$x_{t_n, t_{n+1}}(x, \varepsilon)$ is linear increasing in $x$, and $V^C(x, t_{n+1})$ is concave and increasing in $x$ under the induction assumption of Property 3. Hence, as an increasing concave function of a linear increasing function, $V^C(x_{t_n, t_{n+1}}(x, \varepsilon), t_{n+1})$ is also concave and increasing in $x$. As $x_{t_n, t_{n+1}}(x, \varepsilon) + dt \cdot \sum_{i=n+1}^{N-1} h_{t_i}$, and $V^C(x_{t_n, t_{n+1}}(x, \varepsilon), t_{n+1})$ are concave and increasing in $x$, than $\min\left\{x_{t_n, t_{n+1}}(x, \varepsilon) + dt \cdot \sum_{i=n+1}^{N-1} h_{t_i}, V^C(x_{t_n, t_{n+1}}(x, \varepsilon), t_{n+1})\right\}$ is concave and increasing in $x$, and so as the integral of an increasing concave function. Therefore, we get that $V^C(x, t_n)$ is concave and increasing in $x$. That is, Property 3 is valid for time $t_n$. According to Property 1 and 2, there exist finite prices, $x^H(t_n)$, $x^L(t_n)$ for which any $x \geq x^H(t_n)$, satisfies $x + dt \cdot \sum_{i=n+1}^{N-1} h_{t_i} > V^C(x, t_n)$, and any $x \leq x^L(t_n)$ satisfies $x + dt \cdot \sum_{i=n+1}^{N-1} h_{t_i} < V^C(x, t_n)$. Combined with Property 3 which maintains that $V^C(x, t_n)$ is concave and increasing in $x$, we get that there exists a finite, single value of $x$, defined as the threshold value, $b(t_n)$, for which $b(t_n) + dt \cdot \sum_{i=n}^{N-1} h_{t_i} = V(b(t_n), t_n)$ holds. Figure 5 exemplifies these 3 Properties for $t = 10, T = 15$, $dt = 1, K = 0.5$ and $h_{t_i} = (t_N - t_i) \cdot 10^{-2}$. ∎

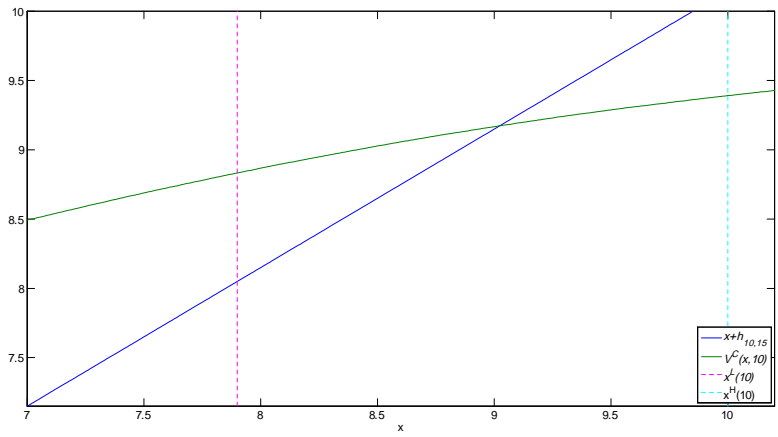

Figure 5: The functions $x + h_{10,15}$ and $V^C(x,t)$ as a function of $x$ at $t = 10, T = 15, \theta = 10, \sigma = 1$, with the bounds $x^L(10), x^H(10)$.